\documentclass[10pt,twocolumn]{article}
\usepackage[margin=1in]{geometry}

\usepackage{amsmath,amssymb,amsthm}
\usepackage{graphicx,color}
\usepackage[hyphens]{url}
\usepackage{dsfont}
\usepackage{booktabs}
\usepackage[normalem]{ulem}

\usepackage{mathtools}
\usepackage[nameinlink,capitalize]{cleveref}
\usepackage[noend]{algpseudocode}
\usepackage{MnSymbol} 

\usepackage[square,sort,numbers]{natbib}
\usepackage[sort,nocompress,noadjust]{cite}

\usepackage{algorithmicx}
\usepackage[Algorithm,ruled]{algorithm}

\usepackage{caption}
\usepackage{subcaption}

\numberwithin{equation}{section}

\newtheorem{lemma}{Lemma}[section]

\theoremstyle{definition}
\newtheorem{defin}{Definition}

\newcommand{\R}{\mathbb{R}}
\newcommand{\C}{\mathbb{C}}

\newcommand{\Rpn}{\R_{\geq 0}^n}
\newcommand{\Rpnn}{\R_{\geq 0}^{n \times n}}

\newcommand{\Rppnn}{\R_{> 0}^{n \times n}}

\providecommand{\diag}{\operatorname{\mathbb{D}}}

\newcommand{\bone}{\mathbf{1}}

\newcommand*{\Otilde}{\tilde{O}}

\newcommand*{\poly}{\mathrm{poly}}

\newcommand*{\eps}{\varepsilon}

\newcommand{\plusminus}{\raisebox{.2ex}{$\scriptstyle\pm$}}

\renewcommand{\geq}{\geqslant}

\newcommand{\NP}{\textsf{NP}}
\newcommand{\OT}{\textsf{OT}}
\newcommand{\MMC}{\textsf{MMC}}

\DeclareMathOperator*{\smin}{smin}
\DeclareMathOperator*{\Dnn}{\Delta_{n \times n}}
\begin{document}

\vspace{-0.7in}
\title{Flows, Scaling, and Entropy Revisited: \\ A Unified Perspective via Optimizing Joint Distributions}
\author{Jason Altschuler\footnote{The author is with the Laboratory for Information and Decision Systems (LIDS), Massachusetts Institute of Technology, Cambridge, MA, 02139.}}
\date{}
\maketitle

\setlength{\tabcolsep}{10pt} 
\renewcommand{\arraystretch}{1.1} 


In this short expository note, we describe a unified algorithmic perspective on several classical problems which have traditionally been studied in different communities. This perspective views the main characters---the problems of Optimal Transport, Minimum Mean Cycle, Matrix Scaling, and Matrix Balancing---through the same lens of \emph{optimization problems over joint probability distributions $P(x,y)$ with constrained marginals}. While this is how Optimal Transport is typically introduced, this lens is markedly less conventional for the other three problems. This perspective leads to a simple and unified framework spanning problem formulation, algorithm development, and runtime analysis. 

Some fragments of this story are classical---for example, the approach for solving Optimal Transport via the Sinkhorn algorithm for Matrix Scaling dates back at least to the 1960s~\citep{Wil69} and by now is well-known in multiple communities spanning economics, statistics, machine learning, operations research, and scientific computing~\citep{SchZen90,Idel16,PeyCut17,Cut13}.

\par Yet, the story described in this note was only recently developed in full---for example, the use of probabilistic inequalities to prove near-optimal runtimes~\citep{AltPar22bal,AltWeeRig17}, and the parallels between Optimal Transport and Minimum Mean Cycle which provide a framework for applying popular algorithmic techniques for the former problem to the latter~\citep{AltPar22mmc}. These developments all hinge on the perspective of optimizing distributions highlighted in this note.
For some problems, this leads to rigorous guarantees that justify algorithms that have long been used in practice and are nowadays the default algorithms in many numerical software packages (e.g., POT, OTT, OTJulia, GeomLoss, MATLAB, R, Lapack, and Eispack); for other problems, this leads to even faster algorithms in practice and/or theory.

The goal of this note is to explain this story with an emphasis on the unifying connections as summarized in Table~\ref{tab:main}. There are several ways to tell this story. We start by introducing Optimal Transport and Minimum Mean Cycle as graph problems (\S\ref{sec:graph}) before re-formulating them as optimization problems over joint distributions (\S\ref{sec:lp}), which makes clear their connections to Matrix Scaling and Balancing (\S\ref{sec:reg}), and naturally leads to the development (\S\ref{sec:alg}) and analysis (\S\ref{sec:analysis}) of scalable algorithms based on entropic regularization. Throughout this narrative, we keep two threads as equals: the fixed marginal setting (Table~\ref{tab:main}, left) and the symmetric marginal setting (Table~\ref{tab:main}, right). As the parallels between these two threads are nearly exact(!), the mental overhead of two threads is hopefully minimal and outweighed by the pedagogical benefit of unifying these two halves of the story---which have often been studied in separate papers and sometimes even separate communities.

\par The presentation of this article is based on Part I of the author's thesis~\citep{Alt22thesis}, and in the interest of brevity, we refer the reader there for proofs and further references (we make no attempt to be comprehensive here given the short length of this note and the immense literature around each of the four problems, e.g., thousands of papers in the past decade just on Optimal Transport).

\begin{table*}[t] 
	\small
	\centering
	\begin{tabular}{c|cc}
		& Fixed marginals & Symmetric marginals \\ \hline
		Linear program  & \textbf{Optimal Transport}    & \textbf{Minimum Mean Cycle}       \\ 
		Entropic regularization & \textbf{Matrix Scaling}           & \textbf{Matrix Balancing}   \\ 
		Polytope & Transportation polytope & Circulation polytope \\
		Simple algorithm & Sinkhorn algorithm & Osborne algorithm \\
		Algorithm in dual & Block coordinate descent & Entrywise coordinate descent \\
		Per-iteration progress & KL divergence & Hellinger divergence \\		
	\end{tabular}
	\caption{
		Although not traditionally viewed in this way, the four bolded optimization problems can be viewed under the same lens: each optimizes a joint probability distribution $P(x,y)$ with similar constraints and objectives. The purpose of this note is to describe the unifying connections in this table and how they can be exploited to obtain practical algorithms with near-optimal runtimes for all four problems. 
	}
	\label{tab:main}
\end{table*}

\paragraph*{Notation.} We associate the set of probability distibutions on $n$ atoms with the simplex $\Delta_n := \{p \in \Rpn : \sum_i p_i = 1\}$, and the set of joint probability distributions on $n \times n$ atoms with $\Delta_{n \times n} := \{P \in \Rpnn : \sum_{ij} P_{ij} = 1\}$. We write $\bone$ to denote the all-ones vector in $\R^n$, and $G = (V,E,c)$ to denote a graph with vertex set $V$, edge set $E$, and edge weights $c : E \to \R$. One caution: we write $\exp[A]$ with brackets to denote the \emph{entrywise} exponential of a matrix $A$. 


\section{Two classical graph problems}\label{sec:graph}

Here we introduce the first two characters in our story: the problems of Optimal Transport ($\OT$) and Minimum Mean Cycle ($\MMC$). We begin by introducing both in the language of graphs as this helps makes the parallels clearer when we transition to the language of probability afterward.

\subsection{Optimal Transport}

In the language of graphs, the $\OT$ problem is to find a flow on a bipartite graph that routes ``supplies'' from one vertex set to ``demands'' in the other vertex set in a minimum-cost way. For convenience, we renormalize the supply/demand to view them as distributions and abuse notation by denoting both vertex sets by $[n] = \{1, \dots, n\}$.


\begin{defin}[Optimal Transport]\label{def:ot}
	Given a weighted bipartite graph $G = ([n] \cup [n],E,c)$ and distributions $\mu,\nu \in \Delta_{n}$, the $\OT$ problem is
	\begin{align}
		\min_{\substack{f : E \to \R_{\geq 0} \\ \sum_{j \in [n]} f(i,j) = \mu_i, \; \forall i \in [n] \\ \sum_{i \in [n]} f(i,j) = \nu_j, \; \forall j \in [n]}} \sum_{e \in E} f(e) c(e)\,.
		\label{eq:ot-def}
	\end{align}
\end{defin}

$\OT$ dates back to Monge in the 18th century when he asked: what is the least-effort way to move a mound of dirt into a nearby ditch of equal volume~\citep{monge1781memoire}? In operations research, $\OT$ appears in textbook problems where one seeks to, e.g., find the minimum-cost way to ship widgets from factories to stores~\citep{BerTsi97}. Recently, $\OT$ has become central to diverse applications in data science---ranging from machine learning to computer vision to the natural sciences---due to the ability of $\OT$ to compare and morph complex data distributions beyond just dirt mounds and widget allocations. Prototypical examples include data distributions arising from point clouds in statistics, images or 3D meshes in computer graphics, document embeddings in natural language processing, or cell phenotyopes or fMRI brain scans in biology. For details on the many applications of $\OT$, we refer to the recent monograph~\citep{PeyCut17}.

\par A central challenge in all these applications is scalable computation. Indeed, data-driven applications require computing $\OT$ when the number of data points $n$ in each distribution is large. Although $\OT$ is a linear program (LP), it is a very large LP when $n$ is in, say, the many thousands or millions, and it is a longstanding challenge to develop algorithms that can compute $\OT$ (even approximately) in a reasonable amount of time for large $n$. See, e.g., the standard texts~\citep{BerTsi97,Sch03,PeyCut17,AhuMagOrl88} for a discussion of the extensive literature on $\OT$ algorithms which dates back to Jacobi in the 19th century.

\subsection{Minimum Mean Cycle}

\begin{table*}[t]
	\small
	\centering
	\begin{tabular}{c|cc}
		& Fixed marginals & Symmetric marginals \\ \hline
		Graphs & Bipartite Flows & Circulations \\ 
		Matrices & Transportation polytope & Circulation polytope   \\ 
		Distributions & Couplings & Self-couplings \\
	\end{tabular}
	\caption{
		Informal dictionary for translating between graphs, matrices, and distributions. See \S\ref{sec:lp} for details.
	}
	\label{tab:gmd}
\end{table*}

We now turn to the second character in our story. Below, recall that a cycle is a sequence of edges that starts and ends at the same vertex.

\begin{defin}[Minimum-Mean-Cycle]\label{def:mmc}
	Given a weighted directed graph $G = (V,E,c)$, the $\MMC$ problem is
	\begin{align}
		\min_{\text{cycle }\sigma \text{ in }G} \frac{1}{|\sigma|} \;\sum_{e \in \sigma} c(e)\,.
		\label{eq:mmc-def}
	\end{align}
\end{defin}

$\MMC$ is a classical problem in algorithmic graph theory which has received significant attention over the past half century due to its many applications. A canonical textbook example is that, at least in an idealized world, finding arbitrage opportunities on Wall Street is equivalent to solving an $\MMC$ problem~\citep[\S24]{CLRS}. Other classic applications range from periodic optimization (e.g., $\MMC$ is equivalent to finding an optimal policy for a deterministic Markov Decision Process~\citep{ZwiPat96}), to algorithm design (e.g., $\MMC$ provides a tractable alternative for the bottleneck step in the Network Simplex algorithm~\citep{GolTar89}), to control theory (e.g., $\MMC$ characterizes the spectral quantities in Max-Plus algebra~\citep{But10}).

\par Just as for $\OT$, a central challenge for $\MMC$ is large-scale computation. The first polynomial-time algorithm\footnote{It is worth remarking that $\MMC$ is polynomial-time solvable while the seemingly similar problem of finding a cycle with minimum (total) weight is $\NP$-hard~\citep[\S8.6b]{Sch03}.} was based on dynamic programming, due to Karp in 1972~\citep{Karp72}. However, its runtime is $O(n^3)$, and an extensive literature has sought to reduce this cubic runtime in both theory and practice; see e.g., the references within the recent papers~\citep{AltPar22mmc,chen2022maximum,GGTW09}.


\section{Graphs, matrices, and distributions}\label{sec:lp}

As written, the optimization problems~\eqref{eq:ot-def} and~\eqref{eq:mmc-def} defining $\OT$ and $\MMC$ appear quite different. Here we describe reformulations that are strikingly parallel. This hinges on a simple but useful connection between graph flows, non-negative matrices, and joint distributions, as summarized in Table~\ref{tab:gmd}. Below, let $C$ denote the $n \times n$ matrix whose $ij$-th entry is the cost $c(i,j)$ if $(i,j)$ is an edge, and $\infty$ otherwise.

\paragraph*{$\OT$ is linear optimization over joint distributions with fixed marginals.}
Consider a feasible flow for the $\OT$ problem in~\eqref{eq:ot-def}, i.e., a flow $f : E \to \R_{\geq 0}$ routing the supply distribution $\mu$ to the demand distribution $\nu$. This flow is naturally associated with a matrix $P \in \R_{\geq 0}^{n \times n}$ whose $ij$-th entry is the flow $f(i,j)$ on that edge. The netflow constraints on $f$ then simply amount to constraints on the row and column sums of $P$, namely $P1 = \mu$ and $P^T1 = \nu$. Thus, $\OT$ can be re-written as the LP
\begin{align}
	\min_{P \in \Dnn \; : \; P1 = \mu, \; P^T1 = \nu} \langle P, C \rangle\,.
	\label{eq:ot-P}
\end{align} 
This decision set $\{P \in \Dnn : P1 = \mu, \;P^T1 = \nu\}$ is called the transportation polytope and can be equivalently viewed as the space of ``couplings''---a.k.a., joint distributions $P(x,y)$ with first marginal $\mu$ and second marginal $\nu$.

\paragraph*{$\MMC$ is linear optimization over joint distributions with symmetric marginals.} $\MMC$ admits a similar formulation by taking an LP relaxation. Briefly, the idea is to re-write the objective in terms of matrices, as above, and then take the convex hull of the discrete decision set. Specifically, re-write the objective $\frac{1}{|\sigma|} \sum_{e \in \sigma} c(e)$ as $\langle P_{\sigma}, C \rangle$ by associating to a cycle $\sigma$ the $n \times n$ matrix $P_{\sigma}$ with $ij$-th entry $1/|\sigma|$ if the edge $(i,j) \in \sigma$, and zero otherwise. 
By the Circulation Decomposition Theorem (see, e.g.,~\citep[Problem 7.14]{BerTsi97}), the convex hull of the set $\{P_{\sigma} : \sigma \text{ cycle}\}$ of normalized cycles is the set $\{P \in \Dnn : P1 = P^T1\}$ of normalized circulations, and so the LP relaxation of $\MMC$ is
\begin{align}
	\min_{P \in \Dnn \; : \; P1 = P^T1 } \langle P, C \rangle\,,
	\label{eq:mmc-P}
\end{align}
and moreover this LP relaxation is exact. For details see, e.g.,~\citep[Problem 5.47]{AhuMagOrl88}. This decision set $\{P \in \Dnn : P1 = P^T1\}$ can be equivalently viewed as the space of ``self-couplings''---a.k.a., joint distributions $P(x,y)$ with symmetric marginals. 
 

\section{Entropic regularization and matrix pre-conditioning}\label{sec:reg}

Together,~\eqref{eq:ot-P} and~\eqref{eq:mmc-P} put $\OT$ and $\MMC$ on equal footing in that they are both LP over spaces of joint distributions $P(x,y)$ with constrained marginals---fixed marginals for $\OT$, and symmetric marginals for $\MMC$. We now move from problem formulation to algorithm development, continuing in a parallel way.

\par The approach discussed in this note, motivated by the interpretation of $\OT$ and $\MMC$ as optimizing distributions, is to use entropic regularization. Namely, add $-\eta^{-1} H(P)$ to the objectives in~\eqref{eq:ot-P} and~\eqref{eq:mmc-P}, where $H(P) = \sum_{ij} P_{ij} \log P_{ij}$ denotes the Shannon entropy of $P$. See Tables~\ref{tab:fixed} and~\ref{tab:sym}, bottom left. This regularization is convex because the entropy function is concave (in fact, strongly concave by Pinsker's inequality). The regularization parameter $\eta > 0$ has a natural tradeoff: intuitively, smaller $\eta$ makes the regularized problem ``more convex'' and thus easier to optimize, but less accurate for the original problem. 

\par But let's step back. Why use entropic regularization? The modern optimization toolbox has many other convex regularizers. The key benefit of entropic regularization is that it reduces the problems of $\OT$ and $\MMC$ to the problems of Matrix Scaling and Matrix Balancing, respectively. This enables the application of classical algorithms for the latter two problems to the former two problems. Below we introduce these two matrix pre-conditioning problems in \S\ref{ssec:reg-def} and then explain this reduction in \S\ref{ssec:reg-red}.

\begin{table*}[t]
	\small
	\centering
	\begin{tabular}{c|ll}
		\textbf{}                  & \multicolumn{1}{c}{\textbf{Primal}} & \multicolumn{1}{c}{\textbf{Dual}} \\ \hline
		\textbf{Optimal Transport}
		& $\min_{P \in \Delta_{n \times n} \; :\; P1 = \mu, P^T1 = \nu} \;\langle P, C \rangle$                         
		&   $\max_{x,y \in \R^n} \min_{ij} (C_{ij} - x_i - y_j) + \langle \mu,x \rangle + \langle \nu,y \rangle$       \\
		\textbf{Matrix Scaling}   
		& $\min_{P \in \Delta_{n \times n} \; :\;  P1 = \mu, P^T1 = \nu} \;\langle P, C \rangle - \tfrac{1}{\eta} H(P)$                         
		&  $\max_{x,y \in \R^n} \smin_{ij} (C_{ij} - x_i - y_j) + \langle \mu,x \rangle + \langle \nu,y \rangle$            
	\end{tabular}
	\caption{
		Primal/dual LP formulations of $\OT$ (top) and its regularization (bottom). The regularization is entropic in the primal and softmin smoothing in the dual. The regularized problem is a convex formulation of the Matrix Scaling problem for the matrix $K = \exp[-\eta C]$.}
	\label{tab:fixed}
\end{table*}

\begin{table*}[t]
	\small
	\centering
	\begin{tabular}{c|ll}
		\textbf{}                  & \multicolumn{1}{c}{\textbf{Primal}} & \multicolumn{1}{c}{\textbf{Dual}} \\ \hline
		\textbf{Minimum Mean Cycle} & $\min_{P \in \Delta_{n \times n} \; :\; P1 = P^T1} \; \langle P, C \rangle $ &   $\max_{x \in \R^n} \min_{ij} C_{ij} + x_i - x_j$       \\
		\textbf{Matrix Balancing}    &  $\min_{P \in \Delta_{n \times n} \; :\;  P1 = P^T1} \;\langle P, C \rangle - \tfrac{1}{\eta} H(P)$  &   $\max_{x \in \R^n} \smin_{ij} C_{ij} + x_i - x_j$           
	\end{tabular}
	\caption{
		Analog to Table~\ref{tab:fixed} in the setting of symmetric marginals rather than fixed marginals. The story is mirrored, with $\OT$ and Matrix Scaling replaced by $\MMC$ and Matrix Balancing, respectively.}
	\label{tab:sym}
\end{table*}

\subsection{Matrix pre-conditioning}\label{ssec:reg-def}

We now introduce the final two characters in our story: Matrix Scaling and Matrix Balancing. In words, these two problems seek to left- and right-multiply a given matrix $K$ by diagonal matrices in order to satisfy certain marginal constraints---fixed marginals for Matrix Scaling, and symmetric marginals for Matrix Balancing. For the former, the scaling is of the form $XKY$; for the latter, it is a similarity transform $XKX^{-1}$.

\begin{defin}[Matrix Scaling]\label{def:ms}
	Given $K \in \Rppnn$ and $\mu,\nu \in \Delta_n$, find positive diagonal matrices $X,Y$ such that $P = XKY$ has marginals $P1 = \mu$ and $P^T1 = \nu$. 
\end{defin}

\begin{defin}[Matrix Balancing]\label{def:mb}
	Given $K \in \Rppnn$, find a positive diagonal matrix $X$ such that $P = XKX^{-1}$ has symmetric marginals $P\bone = P^T\bone$.
\end{defin}

(Both problems are defined here in a simplified way that suffices for the purposes of this note. See the discussion section for details.)

Matrix Scaling and Matrix Balancing are classical problems in their own right and have been studied in many communities over many decades under many names. See the review~\citep{Idel16} for a historical account. The most famous application of these problems is their use as pre-conditioning subroutines before numerical linear algebra computations~\citep{Sin67,Osborne60,PreTeuVetFla07}. For example, Matrix Balancing is nowadays used by default before eigenvalue decomposition and matrix exponentiation in standard numerical packages such as R, MATLAB, Lapack, and Eispack.

\subsection{Reduction}\label{ssec:reg-red}

As alluded to above, entropic regularization leads to the following reductions. Below, $K = \exp[-\eta C]$ denotes the matrix with entries $K_{ij} = \exp(-\eta C_{ij})$. For simplicity, assume henceforth that $G$ is complete so that $K$ is strictly positive, which ensures existence and uniqueness of the Matrix Scaling/Balancing solutions $P$. The general case is similar but requires combinatorial properties of the sparsity pattern~\citep{EavHofRotSch85,rothblum-schneider}. 

\begin{lemma}[Entropic $\OT$ is Matrix Scaling]\label{lem:eot}
	For any $\eta > 0$, the entropic $\OT$ problem has a unique solution. It is the solution $P=XKY$ to the Matrix Scaling problem for $K = \exp[-\eta C]$. 
\end{lemma}

\begin{lemma}[Entropic $\MMC$ is Matrix Balancing]\label{lem:emmc}
	For any $\eta > 0$, the entropic $\MMC$ problem has a unique solution. It is the the solution $P=XKX^{-1}$ to the Matrix Balancing problem for $K = \exp[-\eta C]$, modulo normalizing this Matrix Balancing solution $P$ by a constant so that the sum of its entries is $1$.
\end{lemma}

Both lemmas are immediate from first-order optimality conditions and are classical facts that have been re-discovered many times; see~\citep{Idel16} for a historical account. There is also an elegant dual interpretation. Briefly, entropic regularization in the primal is equivalent to softmin smoothing in the dual, i.e., replacing $\min_i a_i$ by $\smin_i a_i := -\eta^{-1} \log \sum_{i} \exp(-\eta a_i)$ when writing the dual LP in saddle-point form. See Tables~\ref{tab:fixed} and~\ref{tab:sym}, bottom right. Modulo a simple transformation, the optimal scaling matrices $X,Y$ are in correspondence with the optimal solutions to these dual regularized problems---a fact that will be exploited and explained further below.


\section{Simple scalable algorithms}\label{sec:alg}

Since entropic $\OT$ and entropic $\MMC$ are respectively equivalent to Matrix Scaling and Matrix Balancing (Lemmas~\ref{lem:eot} and~\ref{lem:emmc}), it suffices to solve the latter two problems. For both, there is a simple algorithm that has long been the practitioner's algorithm of choice. For Matrix Scaling, this is the Sinkhorn algorithm; for Matrix Balancing, this is the Osborne algorithm. These algorithms have been re-invented many times under different names, see the survey~\citep[\S3.1]{Idel16}. Pseudocode is in Algorithms~\ref{alg:sink} and~\ref{alg:osborne}. For shorthand, we write $\diag(v)$ to denote the diagonal matrix with diagonal $v$, we write $r(P) = P\bone$ and $c(P) = P^T\bone$ to denote row and column sums, and we write $\odot$ and $./$ to denote entrywise multiplication and division.

Both algorithms have natural geometric interpretations as alternating projection in the primal and coordinate descent in the dual. We explain both interpretations as they provide complementary insights.

\begin{algorithm}[h]
	\small
	\centering
	\floatname{algorithm}{\small \textsc{Algorithm}}
	\begin{algorithmic}[1]
		\State Initialize $X,Y \gets I$
					 \Comment{No scaling}
		\State Until convergence:
		\State \indent $X \gets X \odot \mathbb{D}(\mu ./ r(XKY))$ \label{line:sink-rows}\Comment{Fix rows}
		\State \indent $Y \gets Y \odot \mathbb{D}(\nu ./ c(XKY))$ \Comment{Fix columns}
	\end{algorithmic}
	\caption{
		Sinkhorn's Algorithm for scaling a matrix $K$ to have marginals $\mu,\nu$. To solve OT to $\plusminus \eps$, run on $K = \exp[-\eta C]$ where $\eta \approx \eps^{-1} \log n$.
	}
	\label{alg:sink}
\end{algorithm}

\begin{algorithm}[h]
	\small
	\centering
	\caption{
		Osborne's Algorithm for balancing a matrix $K$ to have symmetric marginals. To solve $\MMC$ to $\plusminus \eps$, run on $K = \exp[-\eta C]$ where $\eta \approx \eps^{-1} \log n$.
	}
	\begin{algorithmic}[1]
		\State Initialize $X \gets I$ 
				\Comment{No balancing}
		\State Until convergence:
		\State \indent Choose coordinate $i \in [n]$ to fix
		\State \indent $X_{ii} \gets X_{ii} \cdot \sqrt{c_i(XKX^{-1}) /r_i(XKX^{-1})}$ 
	\end{algorithmic}
	\label{alg:osborne}
\end{algorithm}

\begin{figure}[h]
	\centering
	\includegraphics[width=0.4\linewidth]{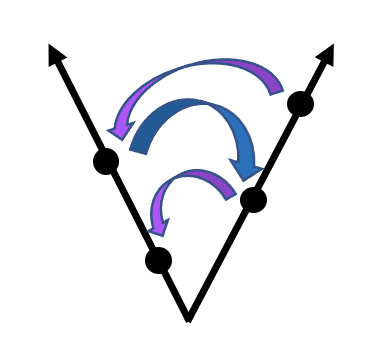}
	\caption{
		In the primal, the Sinkhorn algorithm alternately projects the iterate $P=XKY \in \R^{n \times n}$ onto the two subspaces corresponding to the row/column marginal constraints. The Osborne algorithm is analogous, but with $P=XKX^{-1}$ and $n$ subspaces.}
	\label{fig:primal}
\end{figure}

\paragraph*{Primal interpretation: alternating projection.} The most direct interpretation of the Sinkhorn algorithm is that it alternately projects\footnote{This projection is to the closest point in KL divergence rather than Euclidean distance. In the language of information geometry, this is an I-projection.} the current matrix $P=XKY$ onto either the subspace $\{P \in \Dnn : P1 = \mu \}$ with correct row marginals, or the subspace $\{P \in \Dnn : P^T1 = \nu\}$ with correct column marginals. Note that when it corrects one constraint, it potentially violates the other. Nevertheless, the algorithm converges to the unique solution at the subspaces' intersection, see Figure~\ref{fig:primal} for a cartoon illustration. The Osborne algorithm is analogous, except that it alternately projects $P = XKX^{-1}$ onto $n$ subspaces: the subspaces $\{P \in \Dnn : (P1)_i = (P^T1)_i\}$ defined by equal $i$-th row and column sums, for all $i \in [n]$. Here there is a choice for the order of subspaces to project onto; see~\citep{AltPar22bal} for a detailed discussion of this.

\begin{figure}[h]
	\centering
	\includegraphics[width=0.6\linewidth]{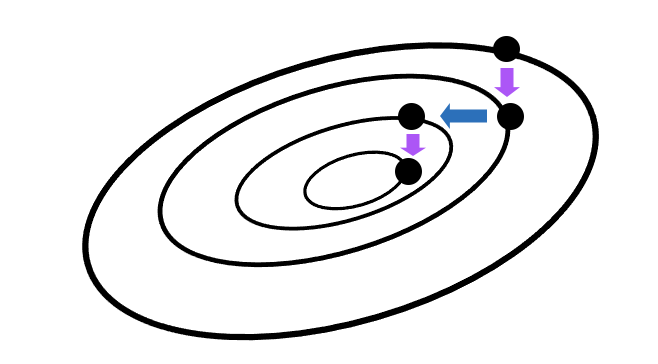}
	\caption{
		In the dual, the Sinkhorn algorithm performs exact block coordinate descent by iteratively updating the scaling matrices $X$ or $Y$ so as to optimally improve the dual objective given that all other entries are fixed. The Osborne algorithm is analogous but updates individual entries rather than blocks.}
	\label{fig:dual}
\end{figure}

\paragraph*{Dual interpretation: coordinate descent.} The Sinkhorn algorithm also admits an appealing dual interpretation. In the dual, entropic $\OT$ has $2n$ variables---the Lagrange multipliers $x,y \in \R^n$ respectively corresponding to the row and column marginal constraints in the primal. See Table~\ref{tab:fixed}, bottom right. Via the transformation $X_{ii} = e^{\eta x_i}$ and $Y_{ii} = e^{ \eta y_i}$, these $2n$ dual variables are in correspondence with the $2n$ diagonal entries of the scaling matrices $X$ and $Y$ that the Sinkhorn algorithm seeks to find. One can verify that the Sinkhorn algorithm's row update (Line~\ref{line:sink-rows} of Algorithm~\ref{alg:sink}) is an exact block coordinate descent step that maximizes the dual regularized objective over all $x \in \R^n$ given that $y$ is fixed. Vice versa for the column update. See Figure~\ref{fig:dual} for a cartoon illustration. The Osborne algorithm is analogous, except that now there are only $n$ dual variables $x \in \R^n$ since there are only $n$ marginal constraints in the primal, see Table~\ref{tab:sym}, bottom right. When the Osborne algorithm updates $X_{ii} = e^{\eta x_i}$ to equate the $i$-th entry of the row and column marginals, this corresponds to an exact coordinate descent step on $x_i$.

\section{Runtime analysis}\label{sec:analysis}

We now turn to convergence analysis. The key challenge is how to measure progress. Indeed, an iteration of the Sinkhorn algorithm corrects either the row or column marginals, but messes up the others. Similarly, an iteration of the Osborne algorithm corrects one row-column pair, but potentially messes up the $n-1$ others. From the cartoons in Figures~\ref{fig:primal} and~\ref{fig:dual}, we might hope to make significant progress in each iteration---but how do we quantify this?

\par There's an entire literature on this question. (Or at least for Matrix Scaling; for Matrix Balancing, things were much less clear until quite recently: even polynomial convergence was unknown for half a century until the breakthrough paper~\citep{OstRabYou16}, let alone near-linear time convergence~\citep{AltPar22bal}.) For example, in the 1980s, Franklin and Lorenz established that Sinkhorn iterations contract in the Hilbert projective metric~\citep{FraLor89}; however, the contraction rate incurs large factors of $n$, which leads to similarly large factors of $n$ in the final runtime. 
Another approach uses auxiliary Lyapunov functions to measure progress, for example the permanent~\citep{LinSamWig98}; however this too incurs extraneous factors of $n$. See the survey~\citep{Idel16}.

\par As it turns out, there is a short, simple, and strikingly parallel analysis approach that leads to runtimes for both the Sinkhorn and Osborne algorithms that scale in the dimension $n$ as $\Otilde(n^2)$~\citep{AltPar22bal,AltWeeRig17}. These are near-optimal\footnote{The ``near'' in ``near-optimal'' refers to the logarithmic factors suppressed by the $\Otilde$ notation.} 
runtimes in $n$ in the sense that in the absence of further structure, it takes $\Theta(n^2)$ time to even read the input for any of the matrix/graph problems discussed in this note, let alone solve them. 

\par At a high level, this analysis hinges on using the regularized dual objective as a Lyapunov function (an idea dating back to the 1990s for Matrix Scaling~\citep{GurYia98}), and observing that each iteration of the Osborne/Sinkhorn algorithm significantly improves this Lyapunov function by an amount related to how violated the marginal constraints are for the current iterate $P$ (``progress lemma''). In its simplest form, the analysis argues that if the current iterate has very violated marginal constraints, then the Lyapunov function improves significantly; and since the Lyapunov function is bounded within a small range (``initialization lemma''), this can only continue for a small number of iterations before we arrive at an iterate with reasonably accurate marginals---and this must be a reasonably accurate solution (``termination lemma'').

\par The progress lemma is itself composed of two steps, both of which leverage the probabilistic perspective highlighted in this note. The first step is a direct calculation which shows that an iteration improves the dual objective by the current imbalance between the marginal distributions as measured in the KL divergence for Sinkhorn, or the Hellinger divergence for Osborne. The second step uses a probabilistic inequality to analyze this imbalance. The point is that since probabilistic inequalities apply for (infinite-dimensional) continuous distributions, they are independent of the dimension $n$. Operationally, this allows switching between the dual (where progress is measured) and the primal (where solutions are desired) without incurring factors of $n$, thereby enabling a final runtime without extraneous factors of $n$. Full details can be found in~\citep{AltWeeRig17,DvuGasKro18} for the Sinkhorn algorithm and~\citep{AltPar22bal} for the Osborne algorithm.


\section{Discussion}

For each of the four problems in this note, much more is known and also many questions remain. Here we briefly mention a few selected topics.

\paragraph*{Algorithm comparisons and the nuances therein.} Each optimization problem in this note has been studied for many decades by many communities, which as already mentioned, has led to the development of many approaches besides the Sinkhorn and Osborne algorithms. Which algorithm is best? 
Currently there is no consensus. We suspect that the true answer is nuanced because different algorithms are often effective for different types of problem instances. For example, specialized combinatorial solvers blow competitors out of the water for small-to-medium problem sizes by easily achieving high accuracy solutions~\citep{dong2020study}, whereas Sinkhorn and Osborne are the default algorithms in most numerical software packages 
for larger problems where moderate accuracy is acceptable.
Both parameter regimes are important, but typically for different application domains. 
For example, high precision may be relevant for safety-critical or scientific computing applications. Whereas the latter regime is typically relevant for modern data-science applications of $\OT$, since there is no need to solve beyond the inherent modeling error ($\OT$ is usually just a proxy loss in machine learning applications) and discretization error ($\mu,\nu$ are often thought of as samples from underlying true distributions). This latter regime is also relevant for pre-conditioning matrices before eigenvalue computation, since Matrix Scaling/Balancing optimize objectives that are just proxies for fast convergence of downstream eigenvalue algorithms~\citep{Osborne60,PreTeuVetFla07}.

\paragraph*{Theory vs practice.}
Comparisons between algorithms are further muddled by discrepancies between theory and practice. Sometimes algorithms are more effective in practice than our best theoretical guarantees suggest---is this because current analysis techniques are lacking, or because the input is an easy benchmark, or because of practical considerations not captured by a runtime theorem (e.g., the Sinkhorn algorithm is ``embarrasingly parallelizable'' and interfaces well with modern GPU hardware~\citep{Cut13}). On the flipside, some algorithms with incredibly fast theoretical runtimes are less practical due to large constant or polylogarithmic factors hidden in the $\Otilde$ runtime. At least for now, this 
includes the elegant line of work (e.g.,~\citep{SpiTen04SDD,CohMadTsiVla17,ZhuLiOliWig17,brand2020bipartitefocs}) 
based on the Laplacian paradigm, which very recently culminated in the incredible theoretical breakthrough~\citep{chen2022maximum} that solves the more general problem of minimum cost flow in time that is almost-linear in the input sparsity and polylogarithmic in $1/\eps$. I am excited to see to what extent theory and practice are bridged in the upcoming years.

\paragraph*{Exploiting structure.} All algorithms discussed so far work for generic inputs. This robustness comes at an unavoidable $\Omega(n^2)$ cost in runtime/memory from just reading/storing the input. This precludes scaling beyond $n$ in the several tens of thousands, say, on a laptop. For larger problems, it is essential to exploit ``stucture'' in the input. What stucture? This is a challenging question in itself because it is tied to the applications and pipelines relevant to practice. For $\OT$, typically $\mu,\nu$ are distributions over $\R^d$ and the cost $C$ is given by pairwise distances, raised to some power $p \in [1,\infty)$. This is the $p$-Wasserstein distance, which plays an analogous role to the $\ell_p$ distance~\citep{Vil03}. Then the $n \times n$ matrix $C$ is implicit from the $2n$ points in $\mu,\nu$. In low dimension $d$, this at least enables reading the input in $O(n)$ time---can $\OT$ also be solved in $O(n)$ time? A beautiful line of work in the computational geometry community has worked towards this goal, a recent breakthrough being $n \cdot (\eps^{-1} \log n)^{O(d)}$ runtimes for $(1\plusminus \eps)$ multiplicatively approximating $\OT$ for $p=1$~\citep{raghvendra2020near,agarwal2022deterministic}. We refer to those papers for a detailed account of this literature and the elegant ideas therein.
Low-dimensional structure can also be exploited by the Sinkhorn algorithm. The basic idea is that the $n^2$ runtime arises only through multiplying the $n \times n$ kernel matrix $K = \exp[-\eta C]$ by a vector---a well-studied task in scientific computing related to Fast Multipole Methods~\citep{beatson1997short}---and this can be done efficiently if the distributions lie on low-dimensional grids~\citep{solomon2015convolutional}, geometric domains arising in computer graphics~\citep{solomon2015convolutional}, manifolds~\citep{AltBacRud19}, or algebraic varieties~\citep{AltPar21ka}. Preliminary numerics suggest that in these structured geometric settings, Sinkhorn can scale to millions of data points $n$ while maintaining reasonable accuracy~\citep{AltBacRud19}. Many questions remain open in this vein, for example graceful performance degradation in the dimension $d$, instance-optimality, and average-case complexity for ``real-world inputs''.

\paragraph*{Optimizing joint distibutions with many constrained marginals.} This note focuses on optimizing joint distributions $P(x,y)$ with $k=2$ constrained marginals, and it is natural to ask about $k \geq 2$. These problems are called Multimarginal $\OT$ in the case of fixed marginals and linear objectives, and arise in diverse applications in fluid dynamics~\citep{Bre08}, barycentric averaging~\citep{AguCar11}, graphical models~\citep{h20gm}, distributionally robust optimization~\citep{Nat18dro}, and much more; see the monographs~\citep{PeyCut17,Pas15,Nat21}. The setting of symmetric marginals arises in quantum chemistry via Density Functional Theory~\citep{cotar2013density,buttazzo2012optimal}. A central challenge in all these problems is that in general, it is intractable even to store a $k$-variate probability distribution, let alone solve for the optimal one. Indeed, a $k$-variate joint distribution in which each variable takes $n$ values is in correspondence with a $k$-order tensor that has $n^k$ entries---an astronomical number even for tiny $n,k =20$, say. As such, a near-linear runtime in the size of $P$ (the goal in this note for $k=2$) is effectively useless for large $k$, and it is essential to go beyond this by exploiting structure via implicit representations. See part II of the author's thesis~\citep{Alt22thesis} and the papers upon which it is based~\citep{AltBoi20motalg,AltBoi21baryalg,AltBoi21mothard,AltBoi22baryhard} for a systematic investigation of what structure enables $\poly(n,k)$ time algorithms, and for pointers to the extensive surrounding literature.

\paragraph*{Matrix Balancing in $\ell_p$ norms.} The original papers~\citep{ParRei69,Osborne60} studied Matrix Balancing in the setting: given $K \in \C^{n \times n}$, find diagonal $X$ such that the $i$-th row and column of $XKX^{-1}$ have equal $\ell_p$ norm, for all $i \in [n]$. Definition~\ref{def:mb} is equivalent for any finite $p$ (which suffices for this note) and elucidates the connection to optimizing distributions. See~\citep[\S1.4]{AltPar22bal} for details. For the case $p=\infty$, similar algorithms have been developed and were recently shown to converge in polynomial time in the breakthrough~\citep{SchSin17}. It would be interesting to reconcile the case $p=\infty$ as the analysis techniques there seem different and the connection to optimizing distributions seems unclear.

\paragraph*{Entropic $\OT$.} In this note, entropically regularized $\OT$ was motivated as a means to an end for computing $\OT$ efficiently. However, it has emerged as an interesting quantity in its own right and is now used in lieu of $\OT$ in many applications due to its better statistical and computational properties~\citep{Cut13,genevay2019sample,MenNil19}, as well as its ability to interface with complex deep learning architectures~\citep{PeyCut17}. Understanding these improved properties is an active area of research bridging optimization, statistics, and applications. We are not aware of an analogous study of entropic $\MMC$ and believe this may be interesting.

\paragraph*{Acknowledgments.} I am grateful to Jonathan Niles-Weed, Pablo Parrilo, and Philippe Rigollet for a great many stimulating conversations about these topics. This work was partially supported by NSF Graduate Research Fellowship 1122374 and a TwoSigma PhD Fellowship.

\small
\bibliographystyle{plainnat}
\bibliography{main}{}

\end{document}